\documentclass[conference]{IEEEtran}

\usepackage{amsmath}
\usepackage{adjustbox}
\usepackage{amssymb}
\usepackage{mathtools}
\usepackage[none]{hyphenat}
\usepackage{graphicx}
\usepackage{tikz}
\usepackage{tikz-network}
\usetikzlibrary{positioning}
\usepackage{array}
\usepackage{booktabs}
\usepackage{soul}
\usepackage{multirow}
\usepackage{multicol}
\usepackage{subfigure}
\def\BibTeX{{\rm B\kern-.05em{\sc i\kern-.025em b}\kern-.08em
		T\kern-.1667em\lower.7ex\hbox{E}\kern-.125emX}}

\DeclareMathOperator*{\argmax}{arg\,max}
\DeclareMathOperator*{\argmin}{arg\,min}

\author{
\IEEEauthorblockN{Raashid Altaf}\\
\IEEEauthorblockA{\textit{Dept. of CSE}, \textit{IIIT Delhi}\\
Delhi, India \\}
\and
\IEEEauthorblockN{Pravesh Biyani}\\
\IEEEauthorblockA{\textit{Dept. of ECE}, \textit{IIIT Delhi}\\
Delhi, India \\}
}

\begin{document}
	\title{Stochastic Trip Planning in High Dimensional Public Transit Network}
	\maketitle
	
	\section{Abstract}
 
	\par This paper proposes a  generalised framework for density estimation in large networks with measurable spatiotemporal variance in edge weights. We solve the stochastic shortest path problem for a large network by estimating the density of the edge weights in the network and analytically finding the distribution of a path. In this study, we employ Gaussian Processes to model the edge weights. This approach not only reduces the analytical complexity associated with computing the stochastic shortest path but also yields satisfactory performance. We also provide an online version of the model that yields a 30 times speedup in the algorithm's runtime while retaining equivalent performance. As an application of the model, we design a real-time trip planning system to find the stochastic shortest path between locations in the public transit network of Delhi. Our observations show that different paths have different likelihoods of being the shortest path at any given time in a public transit network. We demonstrate that choosing the stochastic shortest path over a deterministic shortest path leads to savings in travel time of up to 40\%. Thus, our model takes a significant step towards creating a reliable trip planner and increase the confidence of the general public in developing countries to take up public transit as a primary mode of transportation.

	\section{Introduction}
	\par A trip planning system in public transit aims to provide efficient and practical options for users to navigate a public transportation network. A sound trip planning system is essential to the usability of a public transit network. It condenses information about the entire network into a system accessible to anyone without any requirement of knowledge about the routes, services, or other details of the public transit system in a city.
	
	\par Traditional trip planning approaches can be broadly categorised into two types based on the data used: static and real-time. Static trip planning methods use fixed transit schedules to plan a journey. This approach only works well in cases where the public transit system reliably operates on a schedule, e.g. a metro/subway. Public transit modes such as buses - especially in a developing country like India - due to various operational reasons, do not necessarily adhere to schedule. A real-time trip planning system relying on static data in such cases may provide unreliable and sub-optimal results. This necessitates the usage of real-time data in trip planning for public transit.
	\par A real-time transit feed includes dynamic information about a transit network, such as trip updates and vehicle positions through the GPS devices installed in the transit. The arrival time of transit at stops is generally estimated using this information \cite{fabrikant_2019}\cite{charul}. 
	The travel time of a transit mode between any two stops in the network is the difference between their estimated arrival times (ETA). These travel times are fixed for a given set of ETAs and are used by real-time trip planners. However, transit travel times depend on factors like traffic conditions, bunching etc. and are therefore inherently stochastic. Taking estimated but fixed values of travel times for a journey fails to account for the variance of the travel times experienced in reality. 
	Consequently, the journey planning methods, typically versions of shortest path algorithms, end up being deterministic and face the same pitfalls as the trip planning methods using static pre-set schedules.  
	\par In this paper, we design a predictive model to find the probability distribution of the shortest path in a public transit network with stochastic edge weights. The travel times experienced by a bus between two points in the network are modeled as random variables, representing the real-time variations of the network. Finding the shortest path in this network means predicting the nature of the transit network at a future time instance, which can change with time of day and traffic conditions, making the "shortest" path not unique. Instead, we determine the likelihood of a path being the shortest at a given moment.
	\par Through our work, we redefine the stochastic shortest path problem in the context of a public transit network. The 'shortest' path between two points in a network with stochastic edges is defined as having the maximum optimality index\cite{sigal_stochastic_1980, kamburowski_technical_1985}. The optimality index of a path is traditionally defined as the probability of the path being the shortest among all possible paths for a source-destination pair. We redefine the optimality index as a joint function of the probability of a path being the shortest, as well as the variance of the distribution of the path. In the case of two paths having similar optimality indices, the path with lower variance is recommended to the user. 
	
	\par To solve the stochastic shortest path problem, we model the bus-based public transit network as a weighted directed graph with the bus-stops as nodes, and the edges between the nodes representing the routes and services of the transit. 
	The transit network graph is high dimensional with over a hundred thousand edges. Further, the edges of the transit graph are spatially and temporally correlated. There is also a measurable temporal variance of the edges in the transit network graph.
    Thus, to find the density of a path in the network, we need to find the joint conditional probability density of the corresponding sequence of edges in the transit network. Finding an analytical solution to the stochastic shortest path problem in this scenario is non-trivial due to the scale of the transit network.

    \par Due to the nature of the random variables, we model the distributions of the edge-weights as Random Processes. Furthermore, we define the total cost of a path for a source-destination pair in the transit network graph as the sum of the weights of the edges constituting the path. Thus, the distribution of a total cost of a path in the network is the convolution of conditional densities of the weights of the corresponding edges. We use real-world historical transit data for estimating the probability densities and the correlation of the edge-weights in the transit network. This data is noisy and has missing data values, which occur due to issues such as lack of network connectivity at various locations throughout the city. The task of density estimation in a public transit network is thus a challenging problem from both theoretical and practical perspectives.
    \par In this paper, we model the edge-weights as Gaussian Processes. Gaussian Process Regression is well suited for the task of density estimation in a transit network because:
    \begin{enumerate}
        \item Through the historical data, we observe that the marginal and conditional distribution of the edge-weights exhibits a distribution that can be easily modelled through Gaussian Processes.
        \item The sum of Gaussian random variables is also Gaussian. Therefore, the distribution of the total cost of a path for a source-destination pair is also a Gaussian Process whose parameters can be analytically obtained given the distribution of the edge-weights.
        \item Gaussian Process Regression is well equipped to deal with noisy data and handle missing data values.
    \end{enumerate}
    \par We demonstrate that our model works well in an online setting, reducing resource constraints and enabling us to deploy the model for real-world applications with low computational resources. To the best of our knowledge, this is a first attempt towards solving the stochastic shortest path problem for a large public transit using the real-time and real-world data. A successful implementation will drastically improve the accessibility of public transit for commuters. We also use GTFS, a commonly used data format for open data sources. This ensures that the model can easily be implemented for the transit network of any city.
	\par Our major contributions through this paper include: 
	\begin{enumerate}
		\item Re-look at the stochastic shortest path problem for a public transit network using real-time transit data and find the optimal path in the network for a source-destination pair
		\item Model the transit network as a weighted directed graph with random edge-weights and employ Gaussian processes based density estimation of the edge weights using real-world data. 
		\item Demonstrate the performance, specially in an online setting, as well as the scalability of both probability density estimation as well as real-time journey planning algorithms in a real-world scenario of Delhi with more than two thousand routes, six thousand nodes, and a total of over a hundred throusand edges.
	\end{enumerate}
	\par We first define the transit network and the stochastic shortest path for a public transit network in section \ref{sec:problem_def}. In section \ref{sec:model}, we give a mathematical model for the stochastic shortest path problem in a public transit network. We describe the mathematical formulation of a path in a stochastic network followed by the properties of a path having maximum optimality indes at a given time. We follow this by demonstrating the methodology to implement this trip planning system in a trip planning system in an online setting. Section \ref{sec:experiments} describes the structure and properties of the data used for the experiments. We also detail the analysis performed on the data and describe the challenges faced in pre-processing and estimation phases due to the quality and nature of the data available. The observations and the results are presented in section \ref{sec:observations}
	\subsection{Related Works}
	\par The problem of path-finding in a transit network has seen much research in the field of operations research. Researchers commonly use Dijkstra's Algorithm because of its low complexity and simplicity, enabling researchers to modify the algorithm according to their goals \cite{8093444,bozyigit_public_2018, 5223969}. The goal of a path finding algorithm is to find an "optimal" path to get between two points in the network. In the case of a public transit network, optimality is defined as a combination of factors such as path length, number of transfers and ticket prices \cite{7338605}. The time complexity of a path-finding algorithm in public transit network is especially important for it to be of practical use. To achieve this, researchers model the transit network graph in ways that reduces the search space of the algorithm \cite{5223969, hannah, bast_fast_2007}.
	\par The algorithms currently in place in various trip planning systems such as Google Maps \cite{hannah} are deterministic and designed assuming a static nature of the transit network. Although Google has started using real-time public transit data in 2019 to estimate arrival times of buses \cite{fabrikant_2019}, their trip planning algorithm is inherently deterministic; relying on pre-computations that involve static bus schedules and estimation of ETA from the real-time data \cite{eigenwillig_2016}. Real-time data has also been used to estimate travel-time by using statistical models \cite{xin_model_2014}, neural networks \cite{6709989, liu_short-term_2017} and genetic algorithm \cite{kamel_cgomfp_2004}. While the majority of literature is focussed on estimating the travel time of buses for optimal journey planning \cite{6709989, liu_short-term_2017}, some work has also been done on modelling the passengers' travel time by including factors such as waiting time and time taken to walk to a bus stop  \cite{xin_model_2014}. Some researchers solve a vehicle scheduling problem instead to generate an optimal schedule for the vehicles that leads to an optimal journey for the user in a stochastic network \cite{ricard_predicting_2022}
	\par The earliest available works for a stochastic shortest path (SSP) problem aim to find the distribution of the shortest path in a network having randomly distributed edges \cite{frank_shortest_1969}. Further works develop on this idea by laying down a criterion for optimality; where an optimal path is defined as one that maximises the expectation of a utility function. Elliot and Jerzy \cite{sigal_stochastic_1980, kamburowski_technical_1985}define an optimality index, i.e the probability of a path being shorter than all other possible paths and maximise this index. Other works perform pairwise comparisons between all possible paths to determine the shortest path \cite{abi-char_probability-based_2019}. Recent studies also focus on maximising the probabiility of arriving at the destination on time \cite{RePEc:eee:ejores:v:225:y:2013:i:3:p:455-471, Nikolova06stochasticshortest, NIE2009597} or minimising the expected value of a cost function \cite{jaillet1992shortest,waller2001online, RePEc:eee:ejores:v:225:y:2013:i:3:p:455-471, nikolova2006optimal, thomas2007dynamic, peer2007finding}
	\par The stochastic shortest path problem is a Markov Decision Process. Consequently, the optimal path in a transit network has also been modelled as a state-dependent dynamic system where a policy is either a sequence of services \cite{berczi_stochastic_2017} or a sequence of stops \cite{RePEc:eee:ejores:v:225:y:2013:i:3:p:455-471} that is recommended to the traveller to arrive to the destination on time. Due to the high dimensionality of a public transit network, approximations are also introduced to make a tradeoff between accuracy and run-time \cite{demeyer_dynamic_2014} to maximise the probability of reaching the destination on time \cite{RePEc:eee:ejores:v:225:y:2013:i:3:p:455-471}
	\par In our work, we take the approach of using the real-time transit data to model the public transit network to model the public transit network as a Gaussian Process. We define the optimal path as one that maximises the optimality index \cite{sigal_stochastic_1980}. We also demonstrate the use of real-time and static transit data to perform a series of pre-computations that improve the performance for better practical usage of the journey planning system.
		\section{Problem Definition}\label{sec:problem_def}
		\subsection{Network Definition} \label{subsec: Network Def}
		\par We model the transit network as a graph G (V, E) such that V denotes a set $\{v_1, v_2, \ldots v_n\}$ of bus stops and $E = \{e_1, e_2, \ldots, e_m\}$ denotes the edges between any two stops. We also define a set of $|R|$ routes indexed by unique IDs $R\subset N$.
		\par Each route  $r \in R$ can be defined as a sequence of edges $\{e_{i_1}, e_{i_2}, \dots, e_{i_k}\}$ for $e_{i_j} \in E$, $i_j \in \{1, 2, \ldots, m\}$, and $j \in \{1, 2, ..., k\}$.
		
		\subsection{Stochastic Shortest Path}
		\par Let $w_i(t)$ be the random variable describing the weight of the edge $e_i \in E$ at time $t$ for $i \in \{1, 2, \ldots, m\}$. Also, let $p_{\mathbf{W}}(w_1(t), w_2(t), \ldots, w_m(t))$ be the joint distribution of the random vector $\mathbf{W(t)} \triangleq (w_1(t), w_2(t), \ldots, w_m(t))$.
		\par Let $\tau_{i_j}$ is the time taken to arrive at edge $e_{i_j}$ during a trip. Without loss of generality, we set $\tau_{i_1} = 0$, where $\tau_{i_1}$ is the initial time of the trip, i.e. the time at which the query was made by the user. We define a path $\Pi(s, t)$ from source \textit{s} to destination \textit{t}, where $s, t \in V$, as a sequence of edges $\{e_{i_1}, e_{i_2}, \dots, e_{i_l}\}$ whose total path length is:
		\begin{equation} \label{eq:1}
			|\Pi(s, t)| = \sum_{j = 1}^{l}w_{i_j}(\tau_{i_j})
		\end{equation}
		where:
		\begin{equation}\label{eq:2}
			\tau_{i_j} = \tau_{i_j-1} + w_{i_j-1}(\tau_{i_j-1})
		\end{equation}
		\par Assuming there are $p$ paths $\Pi_1(s, t), \Pi_2(s, t), \ldots, \Pi_p(s, t)$  from source \textit{s} to destination \textit{t}, the shortest path can be defined by a random variable M such that:
		\begin{equation} \label{eq:3}
			M = \min_i(|\Pi_i(s, t)|) \qquad i = (1, \ldots p)
		\end{equation}
		
		The CDF of M can be given by:
		\begin{equation} \label{eq:4}
			\begin{split}
				F_M(m) & = P[M\leq l] \\
				& = 1 - \int_l^{\infty}\ldots \int_l^{\infty} p_{\mathbf{|\Pi|}}\Bigl(\Pi_1, \ldots, \Pi_p\Bigr) d_{\Pi_1},\ldots, d_{\Pi_p} 
			\end{split}
		\end{equation}
		Where $p_{\mathbf{|\Pi|}}\Bigl(\Pi_1, \ldots, \Pi_p\Bigr)$ is the distribution of the random vector $|\mathbf{\Pi}| \triangleq (|\Pi_1(s, t)|, |\Pi_2(s, t)| \ldots, |\Pi_p(s, t)|)$.
		\par As demonstrated by equation \ref{eq:1}, the cost of a path is dependent on the distribution of its edge-weights. In order to determine the distribution of the shortest path, it is necessary to calculate the joint distribution of the edge-weights throughout the network.
		
	\subsection{Objective Function}]\label{subsec:obj-func}
	    Given a source-destination pair $(s, t)$, where $s, t \in V$, suppose there are $k$ possible paths $\Pi_1(s, t), \Pi_2(s, t)$ $\ldots, \Pi_k(s, t)$ in the network having total path lengths $|\Pi_1(s, t)|, |\Pi_2(s, t)|, \ldots, |\Pi_k(s, t)|$ respectively. The \textbf{optimality index}, $C_j$ of a path $\Pi_j(s, t)$, is defined as the probability of $\Pi_j(s, t)$ being the shortest among  $\Pi_1(s, t), \Pi_2(s, t)$ $\ldots, \Pi_k(s, t)$ i.e,
	    \begin{align*}
	        C_j &= P\left[\left|\Pi_j\left(s, t\right)\right| < \left|\Pi_i\left(s, t\right)\right|\right]\\
	        & \qquad \forall i\neq j \quad \& \quad i \in \{1, 2, \ldots, k\}
	    \end{align*}
	    \par Traditionally, the shortest path for the given source-destination pair ($s, t$) is defined as a path with the maximum optimality index, i.e a path $\Pi_p(s, t)$ is the shortest path from source $s$ to destination $t$ iff:
	    \begin{align*}
	        p &= \argmax_{j \in \{1, 2, \ldots, k\}} \{C_j\}
	    \end{align*}
	    \par At a particular time instance, multiple paths may have similar optimality indices (without loss of generality, we define "similar optimality indices" as being within $1\%$ of the maximum optimality index). In such cases we need to consider the variance of the respective paths. The travel time of a path with lower variance is less likely to fluctuate during the course of the trip and is thus preferred. We redefine the 'shortest' path as one with a high probability of being the shortest path while having the least variance.
	    \par Thus, if $\Pi_p(s, t)$ is the path with highest optimality index and $\Pi_{i_1}(s, t), \Pi_{i_2}(s, t), \ldots, \Pi_{i_k}(s, t)$ are the paths with 'similar' optimality indices,  $\Pi_{p'}(s, t)$ is the shortest path from source $s$ to destination $t$ iff:
	    \begin{align*}
	        p' &= \argmin_{j \in \{p, i_1, i_2, \ldots, i_k\}}  \{\sigma_p,\sigma_{i_1}, \ldots, \sigma_{i_k}\}
	    \end{align*}
	    where $\sigma_j$ is the variance of the path $\Pi_j(s, t)$.
	    \par Note that any further mention of a path with "highest optimality index" refers to the path $\Pi_{p'}(s, t)$.
		\section{Model Definition}\label{sec:model}
		\par In a real-world transit network, the time it takes a commuter to travel between two points varies throughout the day based on factors such as traffic conditions and bunching. Consequently, the graph model of the transit network has stochastic edge-weights that can be modelled as random processes. 
		\par In the public transit network, a path between two stops can be defined as a sequence of edges connecting the source to the destination. The distribution of the cost of a path is dependent on the distribution of the edge-weights along the path. Establishing the joint distribution of these edge-weights constitutes a challenging task. 
		\par Due to the nature of flow of traffic in a transit network, we cannot consider the edge-weights to be independent. In fact, the edge-weights have a varying degree of correlation between them depending on their relative geographical locations in the network. Furthermore, correlation between two edge-weights also varies according to the time of day and may have different values based on different time instances throughout the day. This increases the size of the search space needed to model the transit network.
		\par The bus network of Delhi has 6747 stops with over 7000 buses plying on 2000 routes. Considering the size of this network, in addition to the nature of correlation between edge-weights as described above, the problem of modelling the public transit network is computationally very expensive.
		\par We begin by estimating the density of the edge-weights in the network. We experimentally show that a Gaussian Process model works best for this purpose. The distribution of the cost of a particular path is then the convolution of the edge-weights that make up the path, and is also a Gaussian Process. This is followed by covariance estimation required to model the distribution of the path. The shortest path is then analytically calculated from the distributions of all possible paths from the source to the destination.
		\par Gaussian Processes are notoriously expensive to train, scaling with a complexity of $\mathcal{O}(n^3)$ for $n$ observations. Further, as the size of the observations increase, the posterior predictions get slower. To combat this, we demonstrate that the above-mentioned transit model can be easily adapted to an online model for better performance in real-world applications.
		\subsection{Estimating Edge-Weight Density}
        \par We model the edge-weight densities through the observations we generate from the historical data. We assume that the edge weights follow Gaussian Processes and test this assumption with statistical and visual methods. We select 1000 random edges and measure their weights $e(t)$ for $\{0 \leq t \leq 24\}$  where $t$ is the time of day in one-hour bins.
        \par The edge-weight modelling as Gaussian Processes is motivated by the data characteristics and the statistical evidence. We use histograms, Q-Q and P-P plots to visually compare the edge-weight samples with the standard normal distribution. We also apply the Kolgomorov-Smirnov test and calculate the KL Divergence to support this comparison. We provide the details of the methods and results in Appendix \ref{gaussian-exp}.
		\par As the most interesting properties of a Gaussian Process are a result of its covariance function, we use a simple mean function in our Gaussian Process model to reduce the complexity of our estimation process. The mean function returns the mean of observations available to the model.
		\par We set the covariance function to be a sum of different known kernel functions. Specifically, we use exponential squared kernel (equation \ref{eq:kernel}) modelled in a noisy environment. Although non parametric estimation of the kernel function might theoretically result in a more accurate covariance function, we show that the chosen kernel function results in satisfactory results for practical use, at low computation complexity. This is an important distinction considering the dimensions of the solution space.
		\begin{equation}\label{eq:kernel}
			k(x, x') = \sigma^2 \exp \Bigl(-\frac{(x - x')^2}{2l^2}\Bigr)
		\end{equation}
		\par Next we tune the parameters of the kernel function according to the observed data. We find the parameters $\theta'$ that maximise the likelihood $p(\mathbf{y}|X, \theta)$ of the edge-weight density conditional of the observed data $X$.
		\[\theta' = \argmax_\theta p(\mathbf{y}|X, \theta)\]
		\par If we denote the mean $\mu_\theta$ and the covariance function $\Sigma_\theta$ as a function of $\theta$ respectively, the marginal likelihood $p(\mathbf{y}|X, \theta)$ is given by:
		\[p(\mathbf{y}|X, \theta) = \frac{1}{\sqrt{(2\pi)^d|\Sigma_\theta|}}\exp\bigg(-\frac{1}{2}(\mathbf{y}-\mu_\theta)^T\Sigma_\theta^{-1}(\mathbf{y}-\mu_\theta)\bigg)\]
		where $d$ is the dimensionality of the marginal and other symbols have their usual meaning.
		\par We can then find the optimal parameters by minimising the negative log likelihood such that:
		\[\theta' = \argmax_\theta p(\mathbf{y}|X, \theta) = \argmin_\theta (-\log p(\mathbf{y}|X, \theta))\]
		\begin{center}
			\begin{figure}[htb]
			    \centering
			    \caption{kernel fitting example \label{kernel_fit}}
				\includegraphics[width=\textwidth/2, height = 2in]{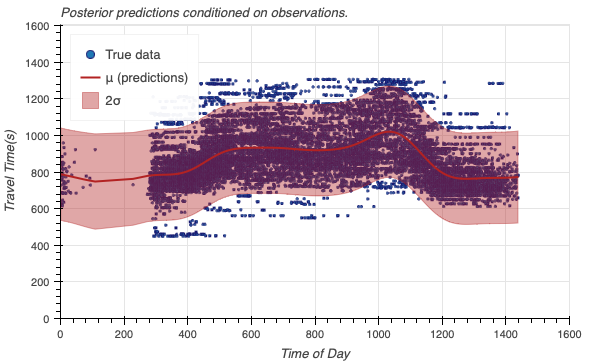}
			\end{figure}
		\end{center}
		\par A gradient based approach is then used to find the optimal parameters (Fig \ref{kernel_fit}). Without loss of generality, we can consider the edge density between an OD pair to be independent of all the routes that pass from the origin to destination. This helps us in designing a much smaller sized network with fewer edges having the same information as all routes passing between two stops can be represented by a single edge. 
		\par Choosing to model edge-weights as Gaussian Processes simplifies the process to obtain the distribution of a path in the transit network. Knowing the parameters of a sequence of edges and the correlation between them, we can easily find the joint distribution of the path.
		\subsection{s-t Path as a Gaussian process}
		\par Note that the edge-weights in the transit network graph are randomly distributed and therefore, the total cost $|\Pi(s, t)|$  of a path $\Pi(s, t) = (e_{i_1}, e_{i_2}, \ldots, e_{i_l})$ is the sum of the weight of the edges (equation \ref{eq:1}):
		 As the edge-weights are modelled as Gaussian Processes, the sum of the edge-weights is also a Gaussian Process $\sim \mathcal{N}(m_{|\Pi|}(t), cov_{|\Pi|}(t, t^*))$ such that:
		\begin{align*} 
			m_{|\Pi|}(t) & = \sum_{j=1}^l m_{w_{i_j}}(t)\\
			cov_{|\Pi|}(t, t^*) &= cov\Bigl(\sum_{j=1}^l w_{i_j}(t), \sum_{j'=1}^l w_{i_{j'}}(t^*)\Bigr)\\
			&= \sum_{j =1}^l \sum_{j' = 1}^l cov(w_{i_j}(t), w_{i_{j'}}(t^*))\\
			&= \sum_{j = 1}^l cov(w_{i_j}(t), w_{i_j}(t^*))\\
			& \qquad + 2\sum_{j<j'} cov(w_{i_j}(t), w_{i_{j'}}(t^*))\\
			\implies cov_{|\Pi|}(t, t^*) &= \sum_{j = 1}^l cov_{w_{i_j}}(t, t^*))\tag{\textbf{I}}\\
			& \qquad + 2\sum_{j<j'} cov(w_{i_j}(t), w_{i_{j'}}(t^*))		
		\end{align*}
		\par Though we can model the density of an edge with the available transit data, it is important to note that the edges in a transit network are not necessarily independent and may depend on other edges in the network spatially as well as temporally. The estimation of covariance between the edge-weights is thus necessary to obtain the distribution of a path in the network (equation I).
		\subsection{Covariance Estimation}
		\par The estimation of the covariance between the edge-weights in a public transit network is a crucial aspect of the probability distribution of the shortest path. A major limitation to using Gaussian Process Regression here is the complexity involved in calculating the covariance matrix. We overcome this by using estimation techniques to obtain the covariance instead of an exact approach. While the variance of each edge-weight can be obtained from the density estimation process, the values of covariance between two different edge-weights are not available a priori. Therefore, in this paper, we estimate these covariance values for every pair of edges for every time instance. This estimation can be performed as a one-time pre-computation. But, to reduce complexity and take advantage of the real-time transit data stream, we employ an online algorithm that updates the measure of covariance between the edges.
		\par To obtain the covariance values, we first calculate the correlation coefficient between the edge-weights and scale them using the variance of the two edge-weights to obtain the covariance (equation \ref{eq:cov}). We calculate the median travel time of an edge $e$ for different hours during the day $t$, over six months of real-time data, Let this be the vector $ETA_e(t)$. We calculate the correlation coefficient for $ETA_{e1}(t)$, and $ETA_{e2}(t')$ for edges $e1$ and $e2$. 
		\begin{equation} \label{eq:cov}
		    corr(x, y) = \frac{cov(x, y)}{\sqrt{var(x)var(y)}}
		\end{equation}
		\par
\begin{figure}[htb]
			\centering
			\caption{Correlation between two consecutive edges w.r.t time of the day \label{corr1}}
			\includegraphics[width=8cm]{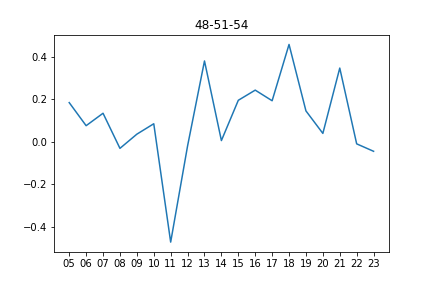}
		\end{figure}

 In Fig \ref{corr1}, the Pearson correlation coefficient is plotted against the time of the day binned by hours. The results are plotted below for two consecutive edges in a route for an instance where $t = t'$
		\par The correlation between two edges in a public transit route refers to the relationship between their travel times. A negative correlation indicates that if one edge is experiencing increased travel times due to high traffic, the subsequent edge on the same route may experience reduced travel times. This can occur when congestion at one location frees up traffic flow for faster speeds further down the route. On the other hand, a positive correlation suggests that both edges are experiencing higher travel times, which is typical during peak hours.
        \par We cannot store the covariance between all the pair of edges for all time due to the size of the resulting dataset. Instead, we create vectors $ETA_e(t)$ for every edge $e$ for all time instances $t$. The value of the covariance is then calculated at run-time as required.
		\par Based on our observations, we have determined that the correlation coefficient changes with time. Additionally, our findings indicate that the edge-weight densities in a transit network are conditionally dependent on each other, in the order of sequence along a route.
		\par It is noteworthy that due to the correlation of the edge-weights, modeling their marginal densities is not sufficient. Hence, in this paper, we also to model the conditional densities of the edge-weights, and determine the shortest path based on the estimated densities.
		\subsection{Shortest Path Estimation}
		\par For every source-destination pair (s, t), where $s, t \in V$, we can obtain a set of possible paths using simple search algorithms such as Depth-First Search(DFS). 
		\par Let the set of paths $\mathbf{\Pi} = \{\Pi_1, \Pi_2, \ldots, \Pi_{m'}\}$\\ for some $s, t \in V$, where
		\[|\Pi_i(t)| \sim \mathcal{GP}(m_i(t), cov(t, t^*)\]
		is a Gaussian process $\forall$ $\Pi_i(t) \in \mathbf{\Pi}$ whose mean and covariance functions are computed apriori as described in the previous sections.
		\par The shortest path will then be the path $\Pi_j$ where:
		\begin{align*}
		    j &= \argmax_j P\bigg[|\Pi_j(s, t)| < \min_{\substack{j\neq i \\ i \in \{1, 2, \ldots k\}}}\bigl\{|\Pi_i(s, t)|\bigr\}\bigg]\\
		    \implies j &= \argmax_jP\Bigl[\bigcap_{i \neq j, i = \{1, ..., m'\}} |\Pi_j| < |\Pi_i|\Bigr]
		\end{align*}
		\par To find the shortest path, we start by finding an initial shortest path from source $s$ to destination $t$. At every possible point of transfer, we find the shortest path from $s'$ to $t$ at $\tau'$, where $s'$ is some stop between $s$ and $t$ and $\tau'$ is the arrival time at $s'$. We suggest a transfer to a different route if we find a better path some time later during the journey. This ensures that we dynamically adjust our results to give optimal results.
		\par
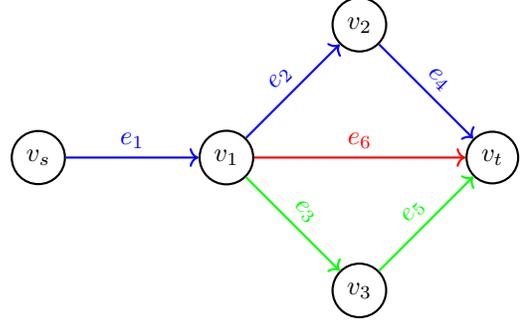
\begin{figure}
			\centering
			\caption{Network Graph} \label{fig:ex_graph}
			\begin{tikzpicture}[node distance={25mm}, thick, main/.style = {draw, circle}]
				\node[main] (1) {$v_s$}; 
				\node[main] (2) [right of=1] {$v_1$}; 	
				\node[main] (3) [above right of=2] {$v_2$};
				\node[main] (4) [below right of=2] {$v_3$};
				\node[main] (5) [above right of=4] {$v_t$};
				\draw[->, blue] (1) -- node[above] {$e_1$} (2);
				\draw[->, red] (2) -- node[above] {$e_6$} (5); 
				\draw[->,blue] (2) -- node[above, sloped] {$e_2$} (3); \draw[->,blue] (3) -- node[above, sloped] {$e_4$} (5);
				\draw[->,green] (2) -- node[above, sloped] {$e_3$} (4); \draw[->,green] (4) -- node[above, sloped] {$e_5$} (5);
			\end{tikzpicture} 
		\end{figure}

 Following example demonstrates this process. Consider the graph in Fig \ref{fig:ex_graph}
		\par Assume that the user starts at $\tau_1 = 0$ and takes route edge $e_1$ with weight density $w_1(\tau_1)$ to reach $v_1$ at time $\lambda_1 = \tau_1 + w_1(\tau_1) = w_1(\tau_1)$. At $v_1$, the user has three options:
		\begin{enumerate}
			\item $\Pi_{l1}(\lambda_1) = w_2(\lambda_1) + w_4(w_2(\lambda_1)) $
			\item $\Pi_{l2}(\lambda_1) = w_3(\lambda_1) + w_5(w_3(\lambda_1)) $
			\item $\Pi_{l3}(\lambda_1) = w_6(\lambda_1) $
		\end{enumerate}
		We choose $\Pi_{li}$ such that:
		\begin{align*}
			i &= \argmax \Bigl[P[\Pi_{l1} < \Pi_{l2} \cap  \Pi_{l1} < \Pi_{l3}], \\
			& \qquad P[\Pi_{l2} < \Pi_{l1} \cap  \Pi_{l2} < \Pi_{l3}], \\
			& \qquad P[\Pi_{l3} < \Pi_{l1} \cap \Pi_{l3} < \Pi_{l2}]\Bigr]
		\end{align*}
		
		Note that all the three options may not be available to the user at any given point. Options 2) and 3) will only be available if a successful transfer takes place from route $l1$ to route $l2$ or $l3$ at stop $v_1$ respectively. 
		
		\par A successful transfer from route $li$ to route $lj$ at stop $v_k$ is said to occur if the arrival time of $li$ at $v_k (\tau^{li}_{v_k})$ is lesser than the arrival time of $lj$ at $v_k (\tau^{lj}_{v_k})$. This transfer can successfully happen with probability 
		\begin{equation} \label{eqn:6}
			P[\tau_{v_k}^{li}\leq \tau_{v_k}^{lj}]
		\end{equation}
		
		\par A problem we face here is that due to the dense nature of the network, we cannot consider the possibility of a transfer at every stop on the initially selected path. Fortunately, we make use of a simple optimisation by only considering stops through which routes can go in multiple directions. More formally, we divide the nodes into two categories: Hub Nodes and Non Hub Nodes:
		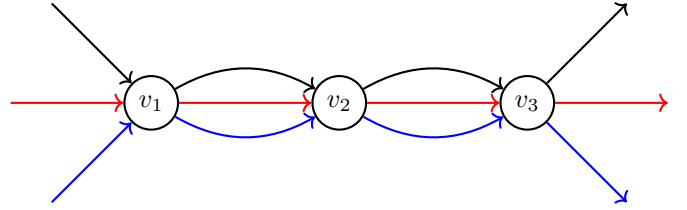
\begin{figure}
			\centering
			\caption{Hub and Non Hub Nodes} \label{fig:hub_nhub}
			\begin{tikzpicture}[node distance={25mm}, thick, main/.style = {draw, circle}]
				\node[main] (1) {$v_1$}; 
				\node[main] (2) [right of=1] {$v_2$}; 	
				\node[main] (3) [right of=2] {$v_3$};
				\node[below left=1.5cm of 1] (4) {};
				\node[above left=1.5cm of 1] (5) {};
				\node[below right=1.5cm of 3] (6) {};
				\node[above right=1.5cm of 3] (7) {};
				\node[left=1.5cm of 1] (8) {};
				\node[right=1.5cm of 3] (9) {};
				\draw[->, blue] (4) -- (1); \draw[->] (5) -- (1);
				\draw[->, blue] (3) -- (6); \draw[->] (3) -- (7);
				\path[->] (1) edge [bend right,blue] (2); \path[->] (1) edge [bend left] (2);
				\path[->] (2) edge [bend right, blue] (3); \path[->] (2) edge [bend left] (3);
				\draw[->, red] (8) -- (1); \draw[->, red] (3) -- (9);
				\draw[->, red] (1) -- (2); \draw[->, red] (2) -- (3);
			\end{tikzpicture}
		\end{figure}
		\subsubsection{Non Hub Nodes:} 
		Node $v_2$ in Fig \ref{fig:hub_nhub} is an example of a non-hub node. These are the nodes in which all the incoming traffic comes from one single direction and goes towards a single direction. As all the traffic moves in a single direction, there's no need to consider any transfer at such stops. 
		\subsubsection{Hub Nodes:}
		Nodes $v_1$ and $v_3$ in Fig \ref{fig:hub_nhub} is an example of a non-hub node. These are the nodes in which either the traffic arrives from multiple directions or departs towards multiple directions, or both.
		
		\par With this information, we formally define the problem of finding a SSP.
		
		\subsection{Stochastic Shortest Path}
		Given s-t paths $\Pi_1(s, t), \Pi_2(s, t)$ $\ldots, \Pi_k(s, t)$ having total path lengths $|\Pi_1(s, t)|, |\Pi_2(s, t)|, \ldots, |\Pi_k(s, t)|$ at a certain point in time $\tau$, respectively. The stochastic shortest path is $\Pi_j(s, t)$ where:
		\[j = \argmax_j P\bigg[|\Pi_j(s, t)| < \min_{\substack{j\neq i \\ i \in \{1, 2, \ldots k\}}}\bigl\{|\Pi_i(s, t)|\bigr\}\bigg]\]
		Now, let
		\begin{eqnarray*}
			F_j&=\bigg\{|\Pi_j(s, t)| < \min_{\substack{j\neq i \\ i \in \{1, 2, \ldots k\}}}\bigl\{|\Pi_i(s, t)|\bigr\}\bigg\}\\
			\implies F_j &= \bigg\{\bigcap_{\substack{j\neq i \\ i \in \{1, 2, \ldots k\}}}|\Pi_j(s, t)| < |\Pi_i(s, t)| \bigg\}
		\end{eqnarray*}
		\par We assume that the densities of all the different paths are pairwise independent. But $\bigl(|\Pi_j(s, t)| < |\Pi_{i_1}(s, t)|\bigr)$ and $\bigl(|\Pi_j(s, t)| < |\Pi_{i_2}(s, t)|\bigr)$ might not necessarily be independent and we do not make any such assumption. Thus to compute $P(F)$ we first find the conditional density $P(F||\Pi_j(s,t)| = \pi)$. We have
		\begin{align*}
			P[F_j||\Pi_j(s, t)| = \pi]&=\\
			&\hspace*{-5mm}P\bigg[\bigcap_{i\neq j}\{|\Pi_j(s, t)| < |\Pi_i(s, t)| \}||\Pi_j(s, t)| = \pi \bigg]\\
			=&P\bigg[\bigcap_{i\neq j}\{\pi < |\Pi_i(s, t)| \}\bigg]\\
			=&\prod_{i\neq j} P[\{|\Pi_i(s, t)| > \pi \}]
		\end{align*}
		\par As we have modelled the path length Gaussian Processes, $|\Pi_i(s, t)| \sim \mathcal{N}(\mu_i(\tau), \sigma_i(\tau))$ is a Gaussian random variable obtained through posterior prediction on the Gaussian Process. We can thus simplify the above equation as
		\begin{align*}
			P[F_j||\Pi_j(s, t)| = \pi]&=\\
			&=\prod_{i\neq j}\bigg[1 - \Phi\bigg(\frac{\pi - \mu_i(\tau)}{\sigma_i(\tau)}\bigg)\bigg]\\
			\implies P[F_j] &=\int P[F_j| |\Pi_j(s, t)| = \pi]f_j(\pi)d\pi\\
			&= \int \prod_{i\neq j}\bigg[1 - \Phi\bigg(\frac{\pi - \mu_i(\tau)}{\sigma_i(\tau)}\bigg)\bigg]f_j(\pi)d\pi\\
			& \text{where }f_j(\pi)\text{ is the pdf of }|\Pi_j(s, t)|
		\end{align*}
		\par There is no closed form solution to the above integral. We thus perform numerical integration $\forall j \in {1, 2, \ldots, k}$ between the 99\% confidence interval of the random variable having the biggest range and select the path with the maximum $P[F_j]$ as the initial shortest path.

        \par Suppose path $\Pi_p\left(s, \right)$ for $p \in {1, 2, \ldots, k}$ is the initial shortest path. As described in section \ref{subsec:obj-func}, we then find the path $\Pi_{p'}\left(s, \right)$ with the highest optimality index. This is the stochastic shortest path between the given source $s$ and the destination $t$.
        
		\subsection{Ranked Shortest Paths}
		\par Although a route might have the highest probability of being the shortest, there are other factor that need to be accounted for to decide which option to suggest to the user. A route can only be the shortest in practice if there is a bus available that traverses that route, in addition to it having the highest probability of being the shortest route. We get this availability information through the estimated arrival time(ETA) for a bus at a particular stop from GPS modules installed on the bus \cite{charul}. 
		\par For any source-destination pair, using the computation provided in the previous section, we can rank all the possible travel options according to the probability of them being the shortest path. We can then choose the route that has the lowest waiting time and the highest probability of being the shortest path.
		\par We can write a path  $\Pi_i(s, t) = (e_{i_1}, e_{i_2}, \ldots, e_{i_m})$ as a sequence of transfers between different routes $e_{i_1}$ to $e_{i_m}$. Let $\eta_{i_{m'}}$ be the earliest ETA of a bus to the head of the edge $e_{i_{m'}},  m' \in {1, 2, \ldots, m}$ and let $e_j(\tau)$ be the density of the edge $e_j$ at time $\tau$  $\forall j \in \{1, 2, \ldots, m\}$
		\par From equation \ref{eq:1} and \ref{eq:2}, we have the total travel time through path $i$, ($tt_i$) given by
		\begin{equation}
			tt_i = \eta_{i_1} + \sum_{j = 1}^m e_{i_j}(\tau_{i_j} + \eta_{i_j})
		\end{equation}
		where $\eta_{i_{m'}} \geq  \sum_{j = 1}^{m'-1} e_{i_j}(\tau_{i_j} + \eta_{i_j})$, i.e. only the buses that arrive at the transfer stop after the user are considered.
		\par From the discussion in the previous section, we can conclude that the path suggested to the user would be the path j such that:
		\begin{equation}\label{eqn:8}
			j = \argmax_j P\big[tt_j < \min_{\substack{i\neq j \\ i \in \{1, 2, \ldots k\}}}\{tt_i\}\big]
		\end{equation}
		\par From here, we can proceed as previous. 
		
		\subsection{Online Learning}
		\par The method proposed above generates results by utilising posterior predictive distributions of the edges. However, this approach presents two challenges. Firstly, it does not allow for the integration of the stream of real-time transit data, as the predictions are solely based on the data used for model training. Secondly, Gaussian Process Regression has a computational complexity of $\mathcal{O}(n^3)$ and a memory complexity of $\mathcal{O}(n^2)$, where $n$ is the size of the training data. In this study, we utilize a large historical transit dataset of 190GB, covering a period of six months, to train our models. This leads to significantly slow predictions that cannot be used in real-world applications. We thus propose an online learning alternative to train Gaussian Process models to counter these two challenges. 
		\par Specifically, we use the Woodbury Idenity with Structured Kernel Interpolation (WISKI)\cite{stanton_kernel_2021} model which combines caching, Woodbury Identiy, and Structured Kernel Interpolation (SKI) to provide constant time (in n) updates while retaining exact inference. Structured Kernel Interpolation (SKI) sparsifies GP through introduction of inducing points. This method proposes an approximation to the kernel matrix $K_{XX} \approx \Tilde{K}_{XX} = WK_{UU}W^\top$, where $U$ is the set of m inducing points, $W \in \mathbb{R}^{n\times m}$ a  sparse cubic interpolation matrix. $W$ consists of $n$ sparse vectors $\mathbf{w}_i \in \mathbb{R}^m$, containing $4^d$ non-zero entries, where $d$ is input dimensions. In SKI, the complexity of adding new datapoints and updating hyperparameters is reduced to $\mathcal{O}(n)$ from $\mathcal{O}(n^3)$ in the native Gaussian Process Regression. This still is not ideal for online learning as the posterior prediction slows down with increase in $n$.
		
		\par WISKI model focuses on reformulating SKI into expressions to get a constant $\mathcal{O}(m^2)$ time and space complexity respectively. Here, Gaussian Process is defined in a regression setting $\mathbf{y} = f(x) + \epsilon$, $f \sim \mathcal{GP}(0, k_\theta(x, x'))$, and $\epsilon \sim \mathcal{N}(0, \sigma^2)$ and $k_\theta(x, x')$ is the kernel function with hyperparameters $\theta$ and $K_{AB}\coloneqq k_\theta(A, B)$ is the covariance between $A$ and $B$. We train the GP hyper-parameters by maximising the marginal log-likelihood using training data $\mathcal{D} = (X, \mathbf{y})$.
		\par The model uses the following equations for to obtain MLL, predictive mean and predictive variance respectively:
		\begin{align*}
		    \log p(\mathbf{y}|X, \theta) &= \frac{1}{2\sigma^2}\big(\mathbf{y^\top}\mathbf{y} - \mathbf{y^\top}WK_{UU}W^\top\mathbf{y} + \\
		    & \mathbf{a^\top}Q^{-1}\mathbf{a}\big) - \frac{1}{2}\left(-\log |Q| + (n-m)\log \sigma^2\right)\\\\
		    \mu_{f|D}\left(\mathbf{x^*}\right) &= \mathbf{w^\top_{x^*}}\left(\sigma^{-2}K_{UU}\left(W^\top\mathbf{y} - L\mathbf{b}\right)\right)\\\\
		    \sigma^2_{f|D}\left(\mathbf{x_i^*}, \mathbf{x_j^*}\right) &= \sigma^2\mathbf{w^\top_{x^*_j}}\left(K_{UU}\left(\mathbf{w_{x^*_j}} - L\mathbf{b'}\right)\right)
		\end{align*}
		\par Here,
		\begin{enumerate}
		\item $LL^\top \approx WW^\top$ is a rank $r$ root decomposition of the matrix $WW^\top$
		\item $Q\coloneqq I + L^\top\sigma^{-2}K_{UU}L$
		\item $\mathbf{b} = Q^{-1}\mathbf{a}$
		\item $\mathbf{a} = L\sigma^{-2}K_{UU}W^\top \mathbf{y}$
		\item $\mathbf{b'} = Q^{-1}L\sigma^{-2}K_{UU}\mathbf{w_{x^*_j}}$
		\item $\mathbf{w}_{t}$ is an interpolation vector for $t$-th data point.
		\end{enumerate}
		\par  For further mathematical explanation and detailed implementation of the WISKI model, we refer the reader to the original paper \cite{stanton_kernel_2021}.
		\par The complexity of computing MLL through this method is $\mathcal{O}(rm + m\log m + jr^2)$ for j steps of conjugate gradients, and $\mathcal{O}(mr^2 + r)$ for conditioning on a new observation. We can see that the total cost depends only on the number of inducing points and the rank of the matrix decomposition. These can be at most $m$, but are typically far less than $m$, which results in constant time updates even as $n$ increases.
		\par A drawback of this method is the memory requirement for higher dimension. According to the authors, if the input data is more than three or four dimensions, the inputs must be projected into a low-dimensional space. But as our data exists in two dimensions, our model doesn't require such projections.
		
		\section{Experiments}\label{sec:experiments}
\begin{figure}[htb]
			\caption{A snapshot of the real-time feed \label{snapshot}}
			\centering
			\includegraphics[width=8cm]{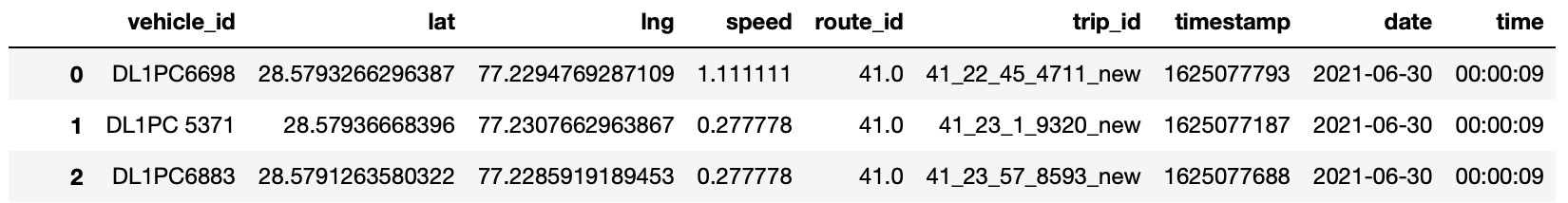}
		\end{figure}

		The data used in the experiments is obtained from the real-time feed of the GPS modules installed on the buses in Delhi. The real-time feed (Fig \ref{snapshot}) contains information about the speed, and location of the bus in addition to identifying information such as the license plate and the route on which the bus is plying. Our APIs fetch this data every 10 seconds throughout the daily service of every bus on the road. We have collected the historical transit data of size 190GB over a period of six months. This section details the techniques employed for the purpose of density estimation of the transit network using this dataset. 
		\subsection{Travel Time}
	    \par As mentioned previously, the time taken by a bus to traverse an edge (hereafter referred to as travel time) is a random variable. To estimate the density function of an edge, we use the historical data to generate samples of travel time for that edge. This information is not directly available to us from the real-time feed, so we compute the approximate travel times of the desired edge for every bus trip.
	    \par If $arr\_time_b(s)$ is the arrival time of bus $b$ at stop $v_s$, then the travel time of a bus $b$ between stop $v_{s_1}$ and $v_{s_2}$, $tt_b(v_{s_1}, v_{s_2})$ can be given by:
	    \[tt_b(v_{s_1}, v_{s_2}) = arr\_time_b(v_{s_1}) - arr\_time_b(v_{s_2})\]
	    \par Note that by taking the difference of arrival times, we implicitly include the time a bus is stationary at a bus stop. 
	    \par As the real-time feed is sampled periodically, it is possible that a bus fails to send any data on its arrival at a stop due to reasons such as network failure or other issues. A reasonable estimate of the arrival time of a bus at a stop, then, is the time at which the distance of the bus from the stop is minimal for a particular trip.
	    \par We can approximate the distance $d_i$ of a bus from the stop at location ($x, y$) at time instance $i$ as the euclidean distance given by
	    \[d_i = \sqrt{(x_i - x)^2 + (y_i - y)^2}\]
	    where ($x_i, y_i$) is the location of the bus at time instance $i$.
	    \par The arrival time of a bus $b$ at a stop $s$ having location ($x, y$) is the time $t$ such that:
	    \[ arr\_time_b(s) = t = \argmin_t d_t\]
	    \par To ensure that the bus is relatively close to a bus stop at the time when the data-point is sampled, we set a threshold on $d_t$. We ignore any trip where $\min d_t > 100$m, i.e if the minima of the distance of a bus from the stop, for a particular trip, is beyond \textbf{100m}, we assume that no data was fetched for that particular trip. The trip is assumed to not provide any significant information, and is discarded. 
	    \par Measuring over a set of 1000 randomly selected edges, we observe that for an edge, approximately 40\% of the total buses passing through a stop send the data within an average minimum distance of 30m $\pm$ 20m from the desired stop. We use the data just from such trips for density estimation purposes.        
		
		\subsection{Simulation Network}
		\par The Delhi public transit system encompasses 6747 bus stops and more than 2000 unique routes that are serviced by 7000 buses daily. As defined in Section \ref{subsec: Network Def}, the network graph of this system results in 6747 nodes and 116316 edges, rendering the computation of edge densities and the application of search algorithms on the entire graph computationally challenging.

        \par Our analysis reveals that the majority of journeys via public transit in Delhi can be completed utilizing a maximum of two buses, or a single transfer. Specifically, a commuter can access an average of 26.74\% of the stops through a single bus, and 99.53\% of the stops with at most one transfer, regardless of their starting point. We use this fact and only consider paths having at most one transfer for every source-destination pair.

        \par Specifically, to optimize the network graph presented in Section \ref{subsec: Network Def}, we establish edges between stops $v$ and $v'$ only if $v'$ is directly accessible from $v$, without any intermediate transfers. As a result, two stops $s$ and $d$ are connected via a single transfer if and only if there exists a stop $v$ such that the edges $s-v$ and $v-d$ are adjacent, with $v$ serving as the transfer point.

        \par This approach reduces the computational complexity of our queries by constructing sub-graphs for all independent origin-destination pairs and estimating the joint density of only two edges, $s-v$ and $v-d$, rather than considering all edges along the path $s-d$. The edge densities in the sub-graphs are then estimated, and the SSP algorithm is run for the desired origin-destination pair.
        
		\par To streamline the presentation, we offer the results for three instances chosen from the 500 pairs tested in certain cases. In particular, we focus on three instances for some metrics and visualizations while providing the complete set of results for other performance measures (Table \ref{instances}). The chosen instances are representative of the algorithm's performance under different scenarios and provide a clear illustration of our findings. The complete set of results is available upon request.
		
        \begin{table}[h]
        \caption{Sample Instances \label{instances}}
        \begin{center}
        \begin{adjustbox}{width = \columnwidth, center }
        \begin{tabular}{lll}
        \toprule
        \textbf{Instance} & \textbf{OD Pair} & \textbf{No of Possible Paths (Transfer Stops)}\\
				\midrule
				1   &  Govind Puri Metro Station to IIT Gate & 18\\
				2   & Govind Puri Metro Station  to ISBT Kashmere Gate Terminal & 14\\
				3   &  Anand Vihar ISBT Terminal to ISBT Kashmere Gate  &  16\\
				\bottomrule
        \end{tabular}
        \end{adjustbox}
        \end{center}
        \end{table}
        
		\begin{figure}[!htb]
		\centering
        \caption{Likelihood of a path through a transfer point being the shortest according to historical transit data \label{p_age_obs}}
        \includegraphics[width=0.9\columnwidth,  height=0.78\textheight]{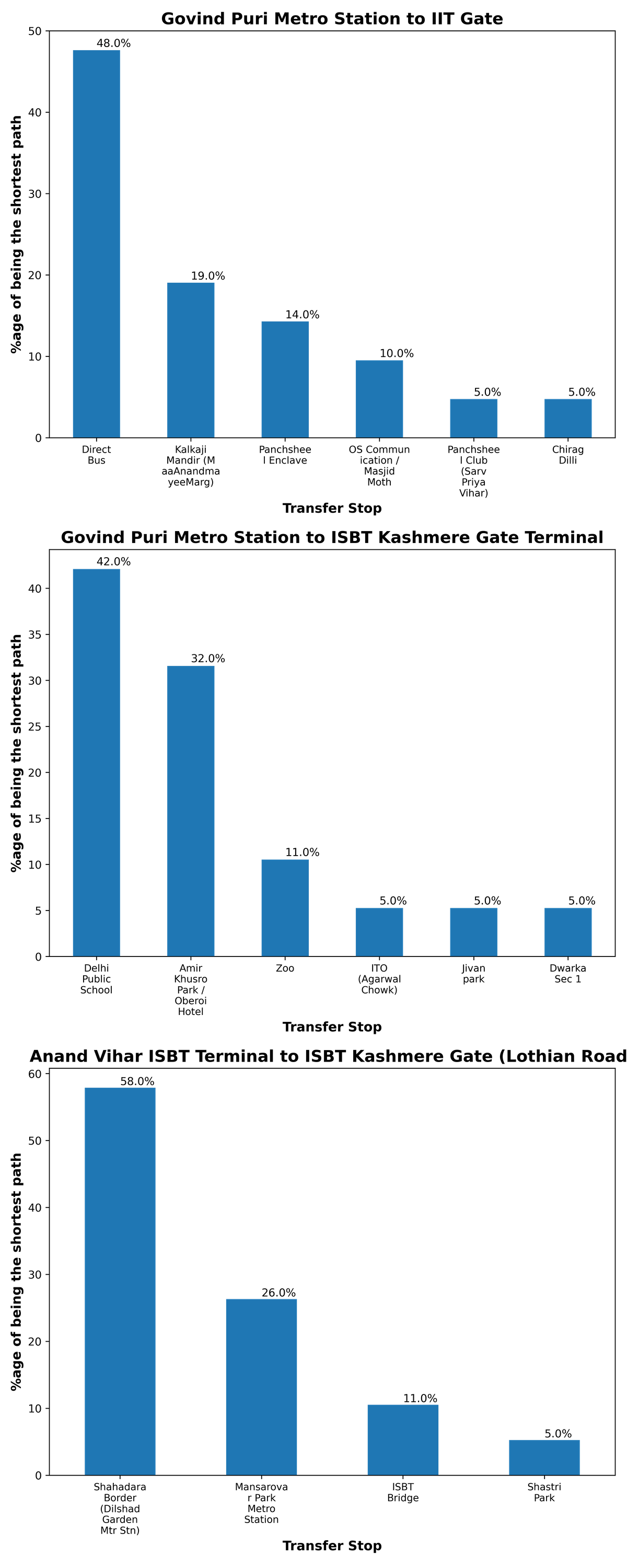}
        \end{figure}
        \begin{figure}[!htb]
        \centering
        \caption{Likelihood of a path through a transfer point being the shortest according to static bus schedules \label{p_age_sch}}
        \includegraphics[width=\columnwidth,  height=0.78\textheight]{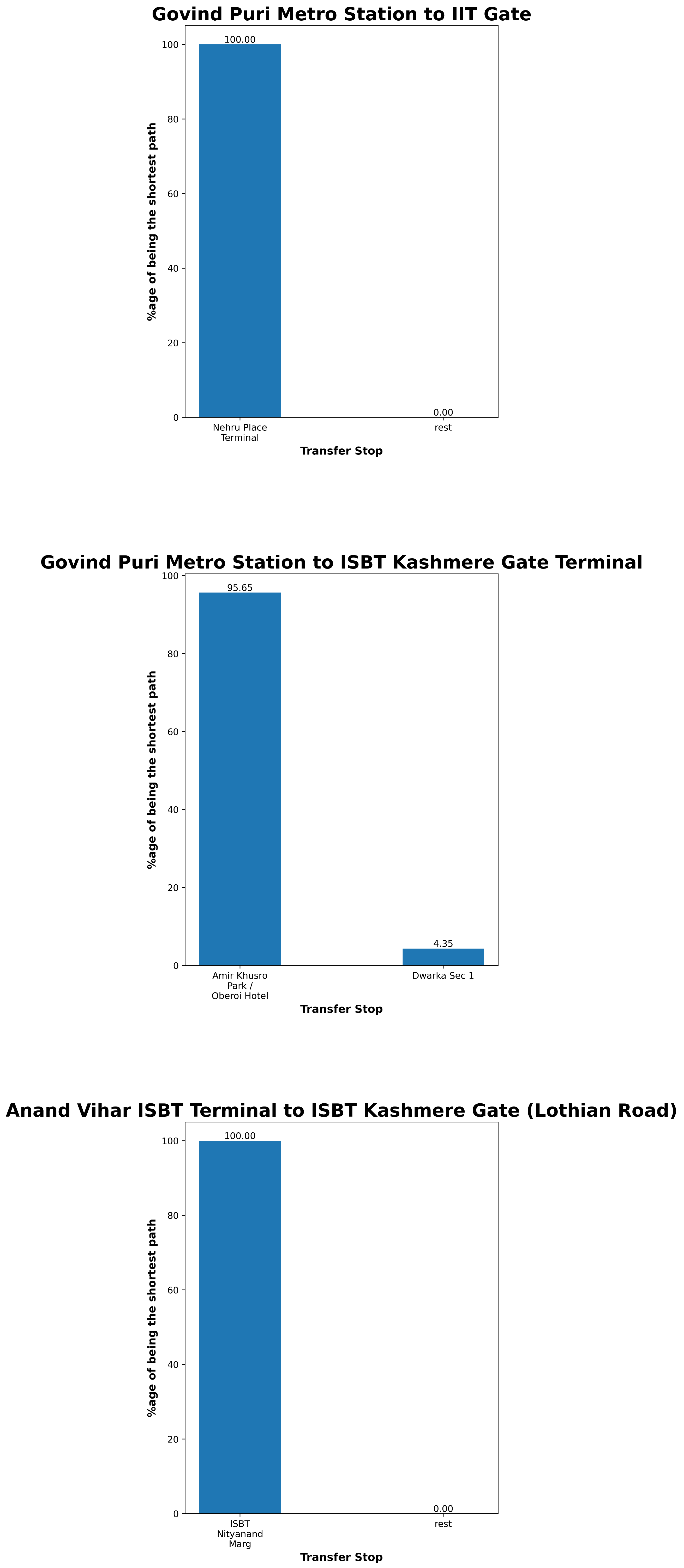}
        \end{figure}
        
		\subsection{Implication of a Stochastic Shortest Path}
		\par This section presents an investigation into the stochastic properties of the shortest path in a public transit network, and the consequential implications for trip planning. We used the historical data to derive travel times for randomly selected origin-destination pairs over a 6-month period at various times throughout the day. Subsequently, we compute the stochastic shortest path for the same times. We uniquely identify a path based on the transfer point between the two legs. The likelihood of a path being the shortest is plotted for selected instances in Fig \ref{p_age_obs}.
		
		\par To compare these results with static trip planning, we also calculate the travel times according to the schedules designed by the transit agencies in Delhi. This is also computed for the same times of day as in the previous case. To maintain parity between the two results, we also include the waiting time between transfers according to the schedules in the total travel time. Fig \ref{p_age_sch} described the results for this experiment. 
		
		\par We observe that the probability of different paths being the shortest vary throughout a 24-hour time period period when considering the historical data while the transit schedule typically results in a single, consistent shortest path throughout the day. We further observe that the shortest paths according to the static bus schedules are often worse than other possible paths in real-life. 
		
		\par Thus, we argue that deterministic calculation of the shortest path is insufficient for trip planning, and a stochastic approach is necessary to obtain accurate results.
        
		\subsection{Experiment Setup}
		\par The models mentioned in this paper are trained using Tensorflow on a laptop equipped with an Apple M1 processor, which features an 8-core CPU, an 8-core GPU, and 16GB of RAM. Python was used along with tensorflow to train the edge-weight density models.
        \subsection{Training Methodologies}
        \par As discussed previously, two methodologies were used to train the models: batch regression and Gaussian Process online learning. Regardless of the approach selected, the process to obtain the SSP remains consistent, with the only difference being the method of training edge-weight densities.
        \par For each OD pair, we perform individual training of each edge in the sub-graph and derive the shortest path results analytically. Consider a path $\Pi_i(s,t)$ comprised of edges $e_1$ and $e_2$. The model was trained to estimate the marginal density of the first edge $p(e_1)$ and the conditional density of the second edge $p(e_2 | e_1)$. The comparison of overall performance of the two methodologies is given in Table \ref{perfcomp}.
        \subsubsection{Batch Regression}
        \par The average training time for all edges in an OD pair sub-graph was 4250.96 seconds, with an average of 140.76 seconds per edge. We reiterate here that the edge-densities were trained on six months of historical transit data. To optimize kernel parameters for each edge, we employed an early stopping method, which involved training the kernel until the loss function converged.
        \par Despite optimization during training, observations show that posterior predictions are significantly slower in the batch regression method (refer Table \ref{perfcomp}). This issue is expected to exacerbate as the model is trained and updated with more transit data. Such extended run-times are not feasible in a practical scenario, even if they provide more accurate results.
    	\begin{table*}[!htbp]
			\centering
			\caption{SSP Observations Sample (Batch Regression) \label{observations}}
			\begin{tabular}{*8c}
				\toprule
				Instance  & No of Edges  &  \multicolumn{2}{c}{Batch Training Time(s)} & \multicolumn{2}{c}{Online Posterior Predictions} & \multicolumn{2}{c}{Batch Posterior Prediction}\\
				\midrule
				&  & Per Edge & Total TIme & Run Time(s) & Probability & Run Time(s) & Probability\\
				1 & 35& 143.79& 5032.76& 0.46& 0.90&215 & 0.93\\
				2 & 28&141.63 & 4451.92&0.40 & 0.76& 37& 0.81\\
				3 & 32 & 140.47 & 4495.20 & 0.41& 0.94& 2087&0.98 \\
				\bottomrule
			\end{tabular}
		\end{table*}
		
        \begin{table}[htbp]
        \centering
        \caption{Performance Comparison Between Online and Batch Regression \label{perfcomp}}
        \begin{tabular}{lll}
        \toprule
        Methodology     & \multicolumn{2}{l}{Average Training} \\ 
        & \multicolumn{2}{l} {Time (seconds)}\\
        \midrule
                        & OD Pair                  & Edge                    \\ 
        Batch Training  & 4250.96                   & 140.76                  \\ 
        Online Training & 150.41                    & 4.29                    \\ 
        \bottomrule
        \end{tabular}
        \end{table}
        \subsubsection{Online Learning}
        \par The Gaussian Process Model was trained online using the WISKI model using the code and methodology presented by Stanton et al. in their publication \cite{stanton_kernel_2021}. To ensure fair comparison with the results of batch regression, the historical data was used in an online setting for the model. 
        \par
\begin{center}
			\begin{figure*}[htb]
					\caption{Mean and variance of the simulation network observed from 0600 to 2300 hours \label{mean_std}}
					\centering
					\includegraphics[height =0.5\textheight]{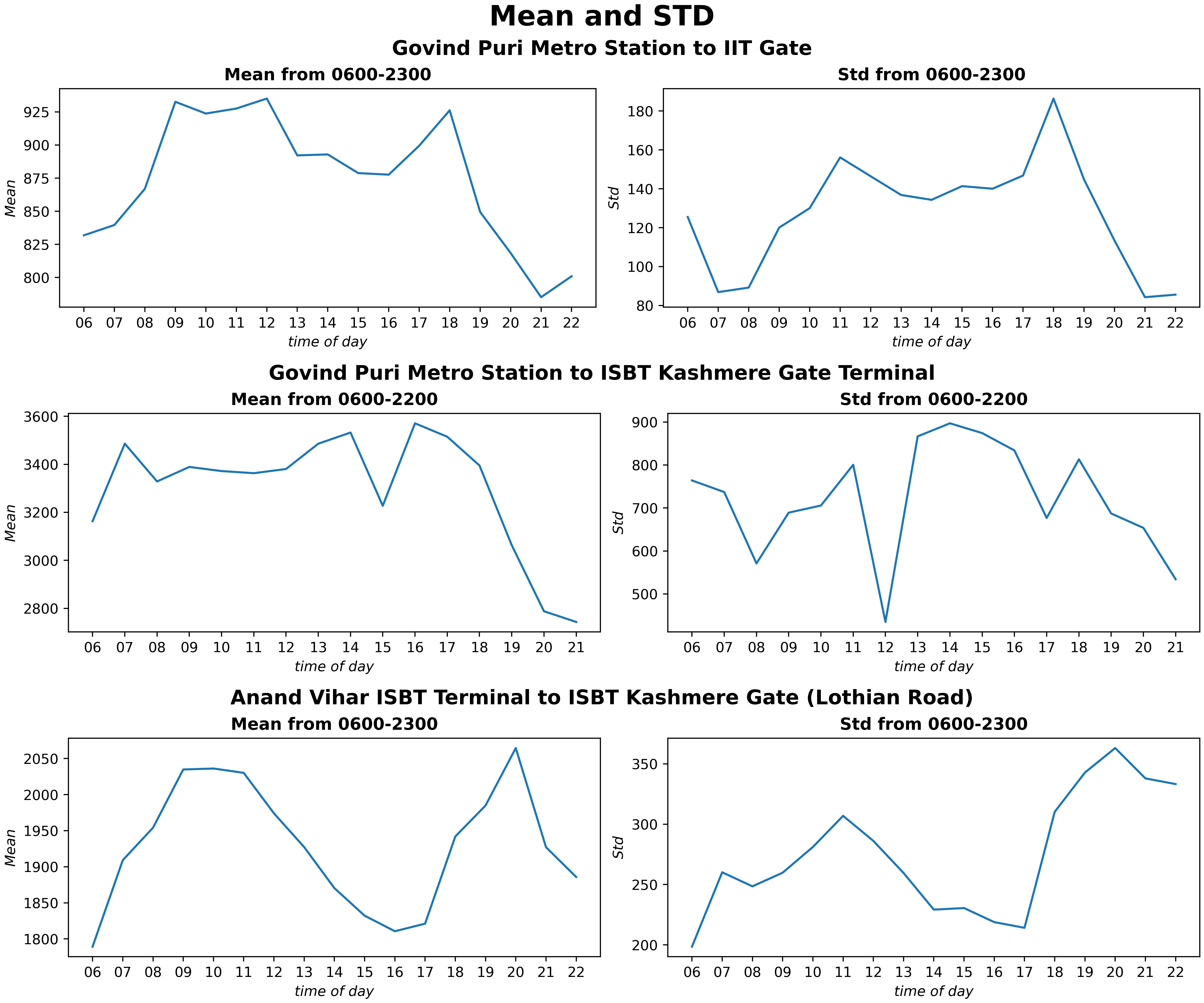}
				\end{figure*}
		\end{center}

 A 95-5\% split was performed on the historical data, with 5\% data being used to train the inital model. The remaining 95\% data was further split into 80-20\% for training and testing, respectively. The results, as shown in Table \ref{training_comp}, demonstrate that the online training approach offers similar performance to batch training while significantly reducing the training time. Specifically, the average training time per edge was 4.3 seconds, representing a 30-fold improvement over batch training. This, in combination with low prediction times (Table \ref{observations}), makes the online training approach suitable for practical applications.
        \begin{table*}[!htbp]
        \centering
		\caption{Online Training Performance Snapshot For a Random Edge \label{training_comp}}
        \begin{tabular}{lllllllllll}
        \toprule
        gp\_loss & batch\_rmse & batch\_nll & online\_rmse & online\_nll & regret & test\_rmse & test\_nll & noise & step\_time & step \\ \midrule
        1.22797 & 1642.19 & 3268.49 & 1669.18 & 3317.66 & 26.9876 & 0.887641 & 1.30686 & 0.664371 & 0.0223811 & 2700 \\ 
        1.22729 & 1703.32 & 3390.46 & 1731.33 & 3440.25 & 28.0081 & 0.888861 & 1.30819 & 0.666535 & 0.0221581 & 2800 \\ 
        1.22999 & 1768.17 & 3521.51 & 1797.04 & 3572.64 & 28.8675 & 0.888199 & 1.30688 & 0.670031 & 0.0216739 & 2900 \\ 
        1.2331  & 1830.75 & 3648.92 & 1860.69 & 3701.26 & 29.9446 & 0.88952  & 1.30817 & 0.673398 & 0.0226879 & 3000 \\ 
        1.22881 & 1885.76 & 3759.48 & 1915.14 & 3812.02 & 29.3783 & 0.887222 & 1.30579 & 0.6687   & 0.0221348 & 3100 \\ \bottomrule
        \end{tabular}
        \end{table*}
		\\\\
		Upon training the kernels, posterior predictions were made on the resulting Gaussian Process Model using the online and batch trained model (as outlined in Equation \ref{eqn:8}). To demonstrate the behaviour of mean and variance of an O-D pair,the results of three instances are displayed in Fig \ref{mean_std}. The peaks observed in the data can be attributed to the rush hour periods. Despite the presence of missing data in the historical data, which results in sharp peaks and dips in the mean and standard deviation plots (as seen in the dip in standard deviation at 12 noon (for instance 2 in Table \ref{instances} of Fig \ref{mean_std}), we believe that with a sufficient amount of time for online learning, the curves will become smoother. The exact cause of this behaviour could not be determined, but with continued training, a more stable and consistent pattern is expected to emerge. 
		\section{Observations}\label{sec:observations}
		\subsection{Discussion of Results}
		\par The performance metric in Table \ref{observations} reflect our observations over all 500 OD pairs we chose for our experiments. We observe that the online and posterior predictions result in different paths having the highest probability of being the shortest path at a given time (Fig \ref{p_age_obs}). Further, we also see that online posterior predictions have a significantly faster run-time than the batch posterior predictions, albeit at the cost of a small drop in the confidence of the results, which is a worthwile exchange.
		\par As we receive a continuous stream of GPS information from the buses, an online model not only leads to low storage use, but also improved performance over the posterior predictive model. 
		
		\par We also draw the readers attention to some peculiar results due to the nature of the data. For instance 2, we observe a lo the low confidence of the posterior predictions. This can be attributed to the high variance of the results for that OD pair. From Fig \ref{mean_std} we can see that the shortest path for instance 2 has a relatively high standard deviation compared to the other two instances. As this implies a more fluctuating travel time, the algorithm has a low confidence in the result. Further, in Fig \ref{p_age_obs} we can see that the difference in likelihood between different options for instance 2 is low as compared to the other 2. This means that there is a higher chance of different paths being the shortest at different times, which is further reflected in the low confidence in the result. This further establishes that our results are within expectations.
        \subsection{Stochastic v/s Static Shortest Path}
        \par To demonstrate the application of stochastic shortest path model in trip planning, we analyse the historical data to evaluate the performance of the proposed stochastic shortest path algorithm against the traditional static schedule approach. Specifically, we compared the actual travel times of the shortest path according to static schedules with those predicted by our model for 20 randomly selected source-destination pairs every hour between 7am and 11pm, for each day in the historical data. The waiting time between the two legs of the journey was also included as a component of the total travel time at each point of transfer. Our findings indicate that the stochastic shortest path algorithm resulted in lower travel times, ranging from 10\% to 40\% lower than the corresponding static shortest path in 96.67\% of the cases. These results provide strong evidence for the potential of the stochastic shortest path approach in improving trip planning.
        
		\par As it stands currently, the online model is capable of generating real-time predictions. In case of systems facing resource constraints, this model may also be used to generate an a-priori ranking of transfer options for all OD paths. The ranks can then be used in addition to the deterministic ETA estimation model \cite{charul} for even lower resource utilization. Exploring that is beyond the scope of this paper.
		
		\section{Conclusion}
		\par In this paper we use a one-of-a-kind historical dataset depicting the traffic pattern of public transit network of Delhi to define the stochastic shortest path problem for a public transit network. Our findings demonstrate that a path in a transit network can be modelled as a Gaussian Process and that the shortest path in the network is stochastic and may change for an origin-destination pair for different times of day. As a result, the likelihood of a path being the shortest is a more accurate measure for trip planning than a deterministic shortest path. 
		\par We model the public transit network in Delhi as a graph, with stops as nodes and bus routes as edges. We utilise the historical dataset, collected by us over a period of six months, consisting of real-time GPS data from the buses in Delhi, to model the edges as independent Gaussian Processes and estimate the correlation between them. This data is noisy and incomplete. To handle these challenges, we employ Gaussian Process Regression for our density estimation process as it is well-suited for this purpose.
		\par Due to the slow posterior predictions in Gaussian Processes, we employ an online learning technique that leads to a drastic increase in training and prediction times while maintaining similar performance. This allows our model to be applicable in real-world use-cases.  
		\par To summarise, the main contribution of our study are the following:
		\begin{enumerate}
		    \item Gathering and using a large real-world transit dataset for modelling transit uncertainty.
		    \item A novel method to model shortest paths in public transit as Gaussian Processes
		    \item Demonstrating that the shortest path in a transit network exhibits a stochastic behaviour.
		    \item Online learning of the Stochastic Shortest Path Problem to achieve milisecond response times. 
		\end{enumerate}
		\par In conclusion, this research highlights the feasibility of using Gaussian Process Regression to tackle the uncertainty present in shortest path problems in public transit networks. With the help of a unique dataset, we have developed a solution that accurately predicts trip plans in real-time. Our findings emphasise the significance of considering transit uncertainty and the necessity for innovative methods to solve such problems. Further studies could focus on scaling up the proposed method for larger transit networks and investigating the possibility of incorporating other sources of uncertainty, such as traffic congestion and road conditions. 		
		\typeout{}
		\bibliographystyle{IEEEtran}
		\bibliography{citations}

\begin{thebibliography}{10}
\providecommand{\url}[1]{#1}
\csname url@samestyle\endcsname
\providecommand{\newblock}{\relax}
\providecommand{\bibinfo}[2]{#2}
\providecommand{\BIBentrySTDinterwordspacing}{\spaceskip=0pt\relax}
\providecommand{\BIBentryALTinterwordstretchfactor}{4}
\providecommand{\BIBentryALTinterwordspacing}{\spaceskip=\fontdimen2\font plus
\BIBentryALTinterwordstretchfactor\fontdimen3\font minus
  \fontdimen4\font\relax}
\providecommand{\BIBforeignlanguage}[2]{{%
\expandafter\ifx\csname l@#1\endcsname\relax
\typeout{** WARNING: IEEEtran.bst: No hyphenation pattern has been}%
\typeout{** loaded for the language `#1'. Using the pattern for}%
\typeout{** the default language instead.}%
\else
\language=\csname l@#1\endcsname
\fi
#2}}
\providecommand{\BIBdecl}{\relax}
\BIBdecl

\bibitem{fabrikant_2019}
\BIBentryALTinterwordspacing
A.~Fabrikant, ``Predicting bus delays with machine learning,'' Jun 2019.
  [Online]. Available:
  \url{https://ai.googleblog.com/2019/06/predicting-bus-delays-with-machine.html}
\BIBentrySTDinterwordspacing

\bibitem{charul}
C.~Paliwal and P.~Biyani, ``To each route its own eta: A generative modeling
  framework for eta prediction,'' in \emph{2019 IEEE Intelligent Transportation
  Systems Conference (ITSC)}, 2019, pp. 3076--3081.

\bibitem{sigal_stochastic_1980}
\BIBentryALTinterwordspacing
C.~E. Sigal, A.~A.~B. Pritsker, and J.~J. Solberg,
  ``\BIBforeignlanguage{en}{The {Stochastic} {Shortest} {Route} {Problem}},''
  \emph{\BIBforeignlanguage{en}{Operations Research}}, vol.~28, no.~5, pp.
  1122--1129, Oct. 1980. [Online]. Available:
  \url{http://pubsonline.informs.org/doi/10.1287/opre.28.5.1122}
\BIBentrySTDinterwordspacing

\bibitem{kamburowski_technical_1985}
\BIBentryALTinterwordspacing
J.~Kamburowski, ``\BIBforeignlanguage{en}{Technical {Note}—{A} {Note} on the
  {Stochastic} {Shortest} {Route} {Problem}},''
  \emph{\BIBforeignlanguage{en}{Operations Research}}, vol.~33, no.~3, pp.
  696--698, Jun. 1985. [Online]. Available:
  \url{http://pubsonline.informs.org/doi/10.1287/opre.33.3.696}
\BIBentrySTDinterwordspacing

\bibitem{8093444}
A.~Bozyiğit, G.~Alankuş, and E.~Nasiboğlu, ``Public transport route
  planning: Modified dijkstra's algorithm,'' in \emph{2017 International
  Conference on Computer Science and Engineering (UBMK)}, 2017, pp. 502--505.

\bibitem{bozyigit_public_2018}
\BIBentryALTinterwordspacing
A.~Bozyigit, E.~Nasiboglu, and S.~Utku, ``Public {Transport} {Route}
  {Recommender} {Regarding} {Multiple} {Factors},'' in \emph{2018 3rd
  {International} {Conference} on {Computer} {Science} and {Engineering}
  ({UBMK})}.\hskip 1em plus 0.5em minus 0.4em\relax Sarajevo: IEEE, Sep. 2018,
  pp. 12--16. [Online]. Available:
  \url{https://ieeexplore.ieee.org/document/8566432/}
\BIBentrySTDinterwordspacing

\bibitem{5223969}
A.~Hedi, Z.~Habbas, and D.~Khadraoui, ``Aco for solving a multimodal transport
  problems using a transfer graph model,'' in \emph{2009 International
  Conference on Computers \& Industrial Engineering}, 2009, pp. 285--290.

\bibitem{7338605}
E.~Nasiboglu, A.~Bozyigit, and Y.~Diker, ``Analysis and evaluation methodology
  for route planning applications in public transportation,'' in \emph{2015 9th
  International Conference on Application of Information and Communication
  Technologies (AICT)}, 2015, pp. 477--481.

\bibitem{hannah}
H.~Bast, E.~Carlsson, A.~Eigenwillig, R.~Geisberger, C.~Harrelson, V.~Raychev,
  and F.~Viger, ``Fast routing in very large public transportation networks
  using transfer patterns,'' vol. 6346, 09 2010, pp. 290--301.

\bibitem{bast_fast_2007}
\BIBentryALTinterwordspacing
H.~Bast, S.~Funke, P.~Sanders, and D.~Schultes, ``\BIBforeignlanguage{en}{Fast
  {Routing} in {Road} {Networks} with {Transit} {Nodes}},''
  \emph{\BIBforeignlanguage{en}{Science}}, vol. 316, no. 5824, pp. 566--566,
  Apr. 2007. [Online]. Available:
  \url{https://www.science.org/doi/10.1126/science.1137521}
\BIBentrySTDinterwordspacing

\bibitem{eigenwillig_2016}
\BIBentryALTinterwordspacing
A.~Eigenwillig, ``An update on fast transit routing with transfer patterns,''
  Mar 2016. [Online]. Available:
  \url{https://ai.googleblog.com/2016/03/an-update-on-fast-transit-routing-with}
\BIBentrySTDinterwordspacing

\bibitem{xin_model_2014}
\BIBentryALTinterwordspacing
G.~Xin and W.~Wang, ``\BIBforeignlanguage{en}{Model {Passengers}’ {Travel}
  {Time} for {Conventional} {Bus} {Stop}},''
  \emph{\BIBforeignlanguage{en}{Journal of Applied Mathematics}}, vol. 2014,
  pp. 1--9, 2014. [Online]. Available:
  \url{https://www.hindawi.com/journals/jam/2014/986546/}
\BIBentrySTDinterwordspacing

\bibitem{6709989}
L.~Deng, Z.~He, and R.~Zhong, ``The bus travel time prediction based on
  bayesian networks,'' in \emph{2013 International Conference on Information
  Technology and Applications}, 2013, pp. 282--285.

\bibitem{liu_short-term_2017}
\BIBentryALTinterwordspacing
Y.~Liu, Y.~Wang, X.~Yang, and L.~Zhang, ``Short-term travel time prediction by
  deep learning: {A} comparison of different {LSTM}-{DNN} models,'' in
  \emph{2017 {IEEE} 20th {International} {Conference} on {Intelligent}
  {Transportation} {Systems} ({ITSC})}.\hskip 1em plus 0.5em minus 0.4em\relax
  Yokohama: IEEE, Oct. 2017, pp. 1--8. [Online]. Available:
  \url{http://ieeexplore.ieee.org/document/8317886/}
\BIBentrySTDinterwordspacing

\bibitem{kamel_cgomfp_2004}
\BIBentryALTinterwordspacing
Z.~Kamel and H.~Slim, ``{CGOMFP} control genetic operators with management of
  the final population to optimize a multimodal transport moving,'' in
  \emph{2004 {IEEE} {International} {Conference} on {Systems}, {Man} and
  {Cybernetics} ({IEEE} {Cat}. {No}.{04CH37583})}, vol.~7.\hskip 1em plus 0.5em
  minus 0.4em\relax The Hague, Netherlands: IEEE, 2004, pp. 6220--6225.
  [Online]. Available: \url{http://ieeexplore.ieee.org/document/1401375/}
\BIBentrySTDinterwordspacing

\bibitem{ricard_predicting_2022}
\BIBentryALTinterwordspacing
L.~Ricard, G.~Desaulniers, A.~Lodi, and L.-M. Rousseau,
  ``\BIBforeignlanguage{en}{Predicting the probability distribution of bus
  travel time to measure the reliability of public transport services},''
  \emph{\BIBforeignlanguage{en}{Transportation Research Part C: Emerging
  Technologies}}, vol. 138, p. 103619, May 2022. [Online]. Available:
  \url{https://linkinghub.elsevier.com/retrieve/pii/S0968090X2200064X}
\BIBentrySTDinterwordspacing

\bibitem{frank_shortest_1969}
\BIBentryALTinterwordspacing
H.~Frank, ``\BIBforeignlanguage{en}{Shortest {Paths} in {Probabilistic}
  {Graphs}},'' \emph{\BIBforeignlanguage{en}{Operations Research}}, vol.~17,
  no.~4, pp. 583--599, Aug. 1969. [Online]. Available:
  \url{http://pubsonline.informs.org/doi/10.1287/opre.17.4.583}
\BIBentrySTDinterwordspacing

\bibitem{abi-char_probability-based_2019}
\BIBentryALTinterwordspacing
P.~E. Abi-Char and A.~Youssef, ``A {Probability}-{Based} {Approach} for
  {Solving} {Shortest} {Path} {Problems} in {Gaussian} {Networks},'' in
  \emph{2019 {IEEE} {Wireless} {Communications} and {Networking} {Conference}
  ({WCNC})}.\hskip 1em plus 0.5em minus 0.4em\relax Marrakesh, Morocco: IEEE,
  Apr. 2019, pp. 1--6. [Online]. Available:
  \url{https://ieeexplore.ieee.org/document/8886027/}
\BIBentrySTDinterwordspacing

\bibitem{RePEc:eee:ejores:v:225:y:2013:i:3:p:455-471}
\BIBentryALTinterwordspacing
L.~Häme and H.~Hakula, ``{Dynamic journeying under uncertainty},''
  \emph{European Journal of Operational Research}, vol. 225, no.~3, pp.
  455--471, 2013. [Online]. Available:
  \url{https://ideas.repec.org/a/eee/ejores/v225y2013i3p455-471.html}
\BIBentrySTDinterwordspacing

\bibitem{Nikolova06stochasticshortest}
E.~Nikolova, J.~A. Kelner, M.~Brand, and M.~Mitzenmacher, ``Stochastic shortest
  paths via quasi-convex maximization,'' in \emph{PROCEEDINGS OF EUROPEAN
  SYMPOSIUM OF ALGORITHMS}.\hskip 1em plus 0.5em minus 0.4em\relax Springer,
  2006, pp. 552--563.

\bibitem{NIE2009597}
\BIBentryALTinterwordspacing
Y.~M. Nie and X.~Wu, ``Shortest path problem considering on-time arrival
  probability,'' \emph{Transportation Research Part B: Methodological},
  vol.~43, no.~6, pp. 597--613, 2009. [Online]. Available:
  \url{https://www.sciencedirect.com/science/article/pii/S0191261509000174}
\BIBentrySTDinterwordspacing

\bibitem{jaillet1992shortest}
P.~Jaillet, ``Shortest path problems with node failures,'' \emph{Networks},
  vol.~22, no.~6, pp. 589--605, 1992.

\bibitem{waller2001online}
S.~T. Waller and A.~K. Ziliaskopoulos, ``On the online shortest path problem,''
  \emph{Science Direct Working Paper}, no. S1574-0358, p.~04, 2001.

\bibitem{nikolova2006optimal}
E.~Nikolova, M.~Brand, and D.~R. Karger, ``Optimal route planning under
  uncertainty.'' in \emph{Icaps}, vol.~6, 2006, pp. 131--141.

\bibitem{thomas2007dynamic}
B.~W. Thomas and C.~C. White~III, ``The dynamic shortest path problem with
  anticipation,'' \emph{European journal of operational research}, vol. 176,
  no.~2, pp. 836--854, 2007.

\bibitem{peer2007finding}
S.~Peer and D.~K. Sharma, ``Finding the shortest path in stochastic networks,''
  \emph{Computers \& Mathematics with Applications}, vol.~53, no.~5, pp.
  729--740, 2007.

\bibitem{berczi_stochastic_2017}
\BIBentryALTinterwordspacing
K.~Bérczi, A.~Jüttner, M.~Laumanns, and J.~Szabó,
  ``\BIBforeignlanguage{en}{Stochastic {Route} {Planning} in {Public}
  {Transport}},'' \emph{\BIBforeignlanguage{en}{Transportation Research
  Procedia}}, vol.~27, pp. 1080--1087, 2017. [Online]. Available:
  \url{https://linkinghub.elsevier.com/retrieve/pii/S2352146517309936}
\BIBentrySTDinterwordspacing

\bibitem{demeyer_dynamic_2014}
\BIBentryALTinterwordspacing
S.~Demeyer, P.~Audenaert, M.~Pickavet, and P.~Demeester,
  ``\BIBforeignlanguage{en}{Dynamic and stochastic routing for multimodal
  transportation systems},'' \emph{\BIBforeignlanguage{en}{IET Intelligent
  Transport Systems}}, vol.~8, no.~2, pp. 112--123, Mar. 2014. [Online].
  Available:
  \url{https://onlinelibrary.wiley.com/doi/10.1049/iet-its.2012.0065}
\BIBentrySTDinterwordspacing

\bibitem{stanton_kernel_2021}
\BIBentryALTinterwordspacing
S.~Stanton, W.~J. Maddox, I.~Delbridge, and A.~G. Wilson,
  ``\BIBforeignlanguage{en}{Kernel {Interpolation} for {Scalable} {Online}
  {Gaussian} {Processes}},'' Mar. 2021, arXiv:2103.01454 [cs, stat]. [Online].
  Available: \url{http://arxiv.org/abs/2103.01454}
\BIBentrySTDinterwordspacing

\end{thebibliography}

    \appendix \label{gaussian-exp}
    \par Consider an edge $e_i \in E$ in a transit network graph $G\left(V, E\right)$. The weight of the edge $w_i$ is the time it takes a bus to traverse the edge. As this time can change throughout the day, we can consider the weight of the edge to be a random process $w(t)$. We assume that the edge-weight density $w_i\left(t\right)$ follows a Gaussian distribution. The following visual and statistical tests performed on 1000 randomly selected edges support this assumption. To ensure proper comparison with the standard normal distribution, the data for edge-weight $w_i\left(t\right) \forall e_i \in E$ was standardised such that it has zero mean and unit variance. 
        \subsection{Visual Tests}
            \begin{enumerate}
                \item Histograms:
                    \par To get a preliminary idea about the properties of the edges we plot the histogram of the edge weights. We select two adjacent edges $e_1(u, v)$ and $e_2(v, w)$ with weights $w_1, w_2$ respectively, having a common vertex $v$. We plot the histograms of $w_1\left(t\right)$ and $w_2\left(t'\right)$ with $t' = t + w_1(t)$. The chosen time $t$ is binned at one-hour time intervals. 
                    \par Upon collecting the data for the required random variables plot the marginal and conditional densities of the edges for different times of the day. Fig \ref{ex1} shows the histogram of the marginal density $p\left(w_1\left(t\right)\right)$ whereas Fig \ref{ex1_cnd} is the histogram of the conditional density $p\left(w_2\left(t'\right)|w_1\left(t\right) = 650 s\right)$ for one such edge pair at time $t = 12$ noon.
        
            		\begin{figure}
            			\caption{Marginal Density $p(w_1)$\label{ex1}}
            			\centering
            			\includegraphics[width=\columnwidth]{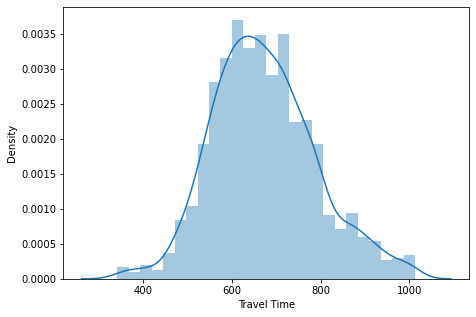}
            		\end{figure}
            		\begin{figure}
            			\caption{Conditional Density $p(w_2|w_1)$\label{ex1_cnd}}
            			\centering
            			\includegraphics[width=\columnwidth]{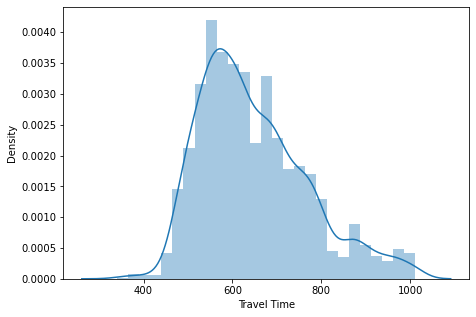}
            		\end{figure}
        
                    \par The edge weights follow a left skewed Gaussian curve with a tail towards the positive end, as shown in the figure. This is because the edge weight model represents travel times, which cannot be negative and thus have outliers only on the positive end of the curve.
                \item Q-Q/P-P Plots:
                    \par A Q-Q plot between two data sets is the plot of the quantiles of the first data set against the quantiles of the second data set. If the two sets come from the same population with the same distribution, the points should fall approximately along a 45-degree reference line plotted along with the data. Similarly, a P-P plot is the plot of the CDFs of the two distributions (empirical and theoretical) against each other. 
        
                    \begin{figure}[htpb]
                        \centering
                        \subfigure[Q-Q Plot]{\includegraphics[width=\linewidth]{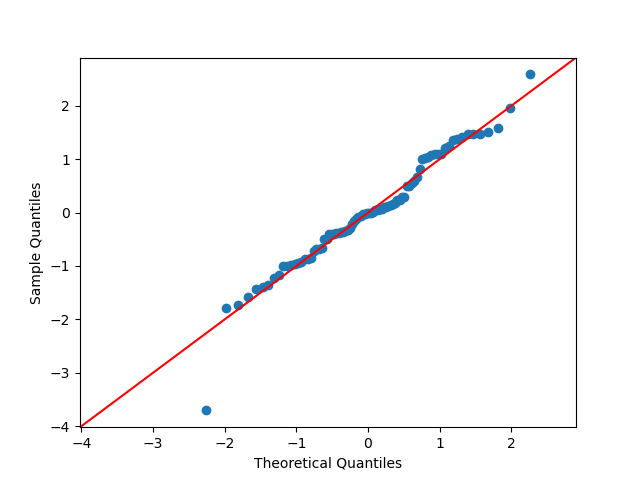}}\\
                        \subfigure[P-P Plot]{\includegraphics[width=\linewidth]{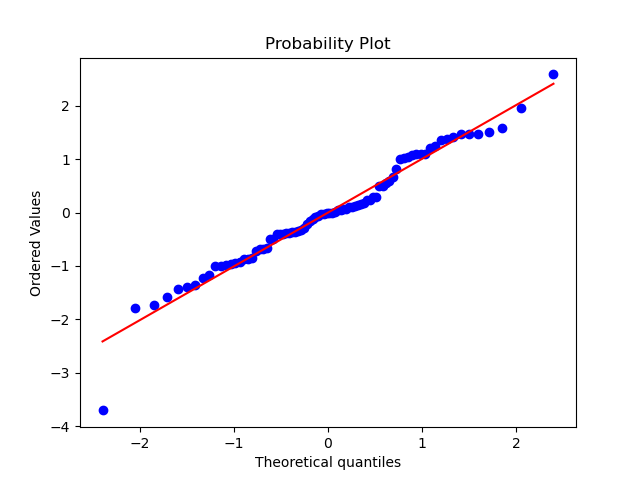}}
                        \caption{Plots for a sample edge with $\mathcal{N}\left(0, 1\right)$ as the theoretical distribution}
                        \label{fig:qpplot}
                    \end{figure}
        
                    \par For both the plots, we plot normalised $w_i(t)$ against a randomly sampled dataset from $\mathcal{N}\left(0, 1\right)$. In Fig. \ref{fig:qpplot}, we have plotted the Q-Q and P-P plots for one instance of edge weight for demonstration. The x-axis in the figures represent the quantiles and the probabilities of the normal distribution respectively whereas the y axes describe the edge-weights. We can see that the datapoints overlap the 45-degree line fairly well except for the deviations towards the ends caused due to outliers. 
        
                    \par The result from the two plots demonstrates that the edge-weights exhibit similar behaviour as is exhibited by a normal demonstration. To further establish this claim, we also perform the following statistical tests.
            \end{enumerate}

        \subsection{Statistical Tests}
            \begin{enumerate}
                \item Kolmogorov-Smirnov Test (KS Test):
                \par The Kolmogorov-Smirnov test is a formal statistical method used to assess the equivalence between continuous one-dimensional probability distributions. It serves the purpose of comparing a given sample to a reference probability distribution. This test aids in determining whether the sample is derived from a population with a particular distribution, which, in our case, is the normal distribution $\mathcal{N}\left(0, 1\right)$.
                \par Formally, we define the KS Test as:
                \begin{itemize}
                    \item $H_0$: The data follow a normal distribution
                    \item $H_a$: The data do not follow a normal distribution
                    \item Test-Statistic: The KS Test statitic is defined as:
                        $$D = \max_{1\leq i\leq N}\left(F\left(Y_i\right) - \frac{i-1}{N}, \frac{i}{N} - F\left(Y_i\right)\right)$$
                        where F is the CDF of the normal distribution
                    \item Significance Level: $\alpha = 0.05$
                    \item Critical Values: We reject the hypothesis regarding the distribution form if the test statistic, D, is greater than the critical value obtained from a table. We use the stats library of Python to perform this test. So, we reject the null hypothesis if the p-value of the test $\Tilde{\alpha} \leq \alpha$
                \end{itemize}
                \par Over the entirety of the dataset of 5000 edge samples, we observe a median $\Tilde{\alpha} = 0.45$, which is far greater than $\alpha$. Thus we cannot reject the null hypothesis.
                \item Kullback-Leibler Divergence:
                \par The Kullback-Leibler (KL) Divergence serves as a metric for quantifying the dissimilarity between two probability distributions. In our case, the objective is to assess the KL divergence between the normalised values of the edge and a standard normal distribution. However, due to the absence of knowledge regarding the underlying distribution from which the edge samples were derived, a simple KL divergence computation solely between the samples and the normal distribution would yield limited information. Furthermore, determining the KL divergence between a distribution and samples obtained from a distribution is not straight-forward.

                \par Therefore, we employ the concept of relative KL divergence to gain insights into the characteristics of the distribution. Specifically, if the KL divergence between two independently sampled standard normal distributions is similar to that between the edge samples $w_i\left(t\right)$ and a sample set from the standard normal distribution, we can reasonably infer the likelihood of the distribution of $w_i\left(t\right)$ being Gaussian.

                \par The following steps were taken to perform this experiment:
                    \begin{enumerate}
                        \item Sample independently from two instances of the standard normal distribution $\mathcal{N}_1(0, 1)$ and $\mathcal{N}_2(0,1)$.
                        \item Calculate the KL Divergence between $\mathcal{N}_1(0, 1)$ and $\mathcal{N}_2(0,1)$ for 10000 iterations to account for the random nature of the data and note the results.
                        \item For every edge $w_i\left(t\right)$, sample randomly from a standard normal distribution $\mathcal{N}'(0, 1)$ a dataset of size equal to $w_i\left(t\right)$ and calculate the KL Divergence between $w_i\left(t\right)$ and $\mathcal{N}'(0, 1)$ and note the results.
                    \end{enumerate}
                \begin{table}[h]
                    \caption{KL Divergence between two randomly sampled standard normal distributions.}\label{tab:kldNorm}
                    \begin{center}
                    \begin{adjustbox}{width = \columnwidth, center }
                    \begin{tabular}{llll}
                    \toprule
                    \textbf{Sample Size} & $\min$ \textbf{KLD} & \textbf{mean KLD} & $\max$ \textbf{KLD} \\ \midrule
                    100         & 0.0006     & 0.036    & 0.432      \\
                    1000        & 0.0007     & 0.007    & 0.046      \\
                    10000       & 0.0004     & 0.001    & 0.002  \\
            				\bottomrule
                    \end{tabular}
                    \end{adjustbox}
                    \end{center}
                \end{table}
                \begin{table}[h]
                     \caption{KL Divergence between edge weight $w_i\left(t\right)$ and $\mathcal{N}\left(0, 1\right)$}\label{tab:kldobs}
                    \begin{center}
                    \begin{adjustbox}{width = \columnwidth, center }
                    \begin{tabular}{llll}
                    \toprule
                    & $\min$ KLD & mean KLD & $\max$ KLD \\ \midrule
                        $KL\left(w_i\left(t\right)||\mathcal{N}\left(0, 1\right)\right)$ & 0.013 & 0.045    & 0.301    \\
            				\bottomrule
                    \end{tabular}
                    \end{adjustbox}
                    \end{center}
                \end{table}
                
                \par Table \ref{tab:kldNorm} and \ref{tab:kldobs} describe the result of this experiment. The size of our observations from $w_i\left(t\right)$ ranged in size form 65 to 105 data points where every data point corresponds to observation from one day. We can see that our observations in Table \ref{tab:kldobs} are consistent with the results in Table \ref{tab:kldNorm} for 100 samples and are within the expected range. 
            \end{enumerate}
            \par To conclude, through the combination of visual and statistical tests, we substantiate our claim that the sample set from $w_i\left(t\right)$ demonstrates similar behaviour to a sample set drawn from a normal distribution. Consequently, modelling the edge-weights as Gaussian Process is a reasonable choice.
            
	\end{document}